%% file: main.tex
\documentclass[runningheads,envcountsame,linenumbers]{llncs}

\usepackage[T1]{fontenc}
\usepackage[utf8]{inputenc}
\usepackage{todonotes}
\usepackage{commath}
\usepackage{amsmath, amsfonts, amssymb, dsfont}
\usepackage{mathtools}
\usepackage{csquotes}
\usepackage{tikz}
\usepackage{hyperref}
\usepackage{cleveref}
\usetikzlibrary{shapes.geometric, fit, positioning}

\usepackage{lineno}

\usepackage{enumitem}
\usepackage{listings}

\input{shortcuts}

\newif\ifpaper
\paperfalse 

\begin{document}

\title{Jumbled Scattered Factors}

\author{Pamela Fleischmann\inst{1}\orcidID{0000-0002-1531-7970} \and
		Annika Huch\inst{1}\orcidID{0009-0005-1145-5806}\and
		Melf Kammholz\inst{1}\and
		Tore Koß \inst{2}\orcidID{0000-0001-6002-1581}}
	\authorrunning{Fleischmann  et al.}
	%
	\institute{Kiel University, Germany,
		\email{\{fpa,ahu\}@informatik.uni-kiel.de, stu209214@mail.uni-kiel.de}\and University of Göttingen, Germany, \email{tore.koss@cs.uni-goettingen.de}}

\maketitle

\input{abstract}

\section{Introduction}
\input{intro}

\section{Preliminaries}
\input{prelims}


\section{The Jumble Index and Jumble Division}\label{sec:indexdivision}
\input{indexdivision}

\section{Jumbled Simon's Equivalence $\jsimkl$}\label{sec:jumblesimon}
\input{jumblesimon}

\section{Jumbled Universality}\label{sec:universality}
\input{universality}
\section{Conclusion}
\input{conclusion}
%
\newpage
\bibliographystyle{plainurl}
\bibliography{refs}

\newpage




\end{document}

%% file: shortcuts.tex
\def\ta{\mathtt{a}}
\def\tb{\mathtt{b}}
\def\tc{\mathtt{c}}
\def\td{\mathtt{d}}

\def\tx{\mathtt{x}}
\def\ty{\mathtt{y}}

\def\nth#1{#1$^{\text{th}}$}

\def\N{\mathbbm{N}}

\def\abs#1{|#1|}

\renewcommand{\epsilon}{\varepsilon}

\renewcommand{\phi}{\varphi}

\def\N{\mathbb{N}} 
\def\Z{\mathbb{Z}} 

\DeclareMathOperator{\al}{alph}
\DeclareMathOperator{\Fact}{Fact}

\usepackage{bbm}

\newcommand{\ScatFact}{\operatorname{ScatFact}}
\newcommand{\JumbledScatFact}{\operatorname{JScatFact}}
\newcommand{\JScatFact}{\operatorname{JScatFact}}
\newcommand{\JSF}{\operatorname{JScatFact}}
\DeclareMathOperator{\SJSF}{\overline{\JScatFact}}
\DeclareMathOperator{\SJScatFact}{\overline{\JScatFact}}

\newcommand{\CSF}{\operatorname{CSF}}
\newcommand{\lcsf}{\operatorname{lcsf}}
\newcommand{\LCSF}{\operatorname{LCSF}}

\newcommand{\ar}{\operatorname{ar}}

\renewcommand{\r}{\operatorname{r}}
\DeclareMathOperator{\letters}{alph}

\newcommand{\jsimkl}{\overset{J}{\sim}_{k,\ell}}
\newcommand{\jsim}{\overset{J}{\sim}}
\newcommand{\jdiv}{\mid_J}

\usepackage{shuffle}

\newcommand\restr[2]{{
  \left.\kern-\nulldelimiterspace 
  #1 
  \vphantom{\big|} 
  \right|_{#2} 
  }}

\newcommand{\cmpl}[1]{#1^\complement}

\newcommand{\newterm}[1]{\textnormal{#1}}

\usepackage{todonotes}

%% file: abstract.tex
\begin{abstract}
In this work, we combine the research on (absent) scattered factors with the one
of jumbled words. For instance, $\mathtt{wolf}$ is an absent scattered factor of
$\mathtt{cauliflower}$ but since $\mathtt{lfow}$, a
jumbled (or abelian) version of $\mathtt{wolf}$, is a scattered factor, $\mathtt{wolf}$ occurs as a jumbled scattered factor in
$\mathtt{cauliflower}$. A \emph{jumbled scattered factor} $u$ of a word $w$ is
constructed by letters of $w$ with the only rule that the number of occurrences
per letter in $u$ is smaller than or equal to the one in $w$. We proceed to
partition and characterise the set of jumbled scattered factors by the number of
jumbled letters and use the latter as a measure. For this new class of words, we
relate the folklore longest common subsequence (scattered factor) to the number
of required jumbles. Further, we investigate the smallest possible number of
jumbles alongside the jumbled scattered factor relation as well as Simon's congruence
from the point of view of jumbled scattered factors and jumbled universality.
\end{abstract}

%% file: intro.tex
In many areas of computer science, bioinformatics and data science, it is
crucial to analyse data. One fundamental concept reoccurring in all these fields
is the concept of \emph{scattered factors} (also known as \emph{scattered subwords} or
\emph{subsequences}), which are derived by selecting non-necessarily consecutive
parts of a string while preserving their relative order. For example, the words
$\mathtt{cute}, \mathtt{outer}$ and $\mathtt{mute}$ are scattered factors of
$\mathtt{computer}$, while $\mathtt{route}$ and $\mathtt{tempo}$ are not since
the letters occur in the wrong order. The study of scattered factors dates
back to the 1970s, when Imre Simon introduced them in the context of
piecewise-testable languages \cite{DBLP:conf/automata/Simon75}. Words and their
scattered factors have been studied in combinatorics on words 
\cite{DBLP:journals/fuin/KoscheKMS22,DBLP:journals/jcss/MateescuSY04},
algorithms and stringology \cite{DBLP:journals/tcs/Baeza-Yates91,DBLP:journals/algorithmica/BannaiIKKP24,DBLP:conf/soda/BringmannK18,DBLP:journals/aam/ForresterM23,DBLP:journals/jacm/Maier78,DBLP:journals/jacm/WagnerF74},
logics \cite{DBLP:journals/lmcs/BaumannGTZ23,DBLP:conf/lics/HalfonSZ17,10.1007/978-3-030-50026-9_21,10.1007/978-3-030-17127-8_20},
as well as language and automata theory \cite{DBLP:conf/isaac/AdamsonFHKMN23,DBLP:journals/ipl/KarandikarKS15,DBLP:journals/lmcs/KarandikarS19,simon1972hierarchies,DBLP:conf/automata/Simon75,zetzsche:LIPIcs.ICALP.2016.123}
(see \cite[Chapter 6, Subwords]{lothaire} for an overview). Since then, scattered
factors played an important role also in a broader variety of fields such as
bioinformatics \cite{DBLP:journals/bib/DoLL21,DBLP:conf/aaai/WangCG20}, formal
software verification \cite{DBLP:conf/lics/HalfonSZ17,zetzsche:LIPIcs.ICALP.2016.123}
and database theory \cite{artikis2017complex,FrochauxK23,Kleest-Meissner22,Kleest-Meissner23,SchmidSIGMOD}.
Further, scattered factors are also used to model corrupted data in the context
of the theoretical problem of reconstruction
\cite{dress2005reconstructing,DBLP:journals/ijfcs/FleischmannLMNR21,DBLP:conf/dlt/Manuch99}.
The investigation of Simon's congruence
\cite{DBLP:conf/automata/Simon75}, where two words are called
\emph{Simon congruent} if they have the same sets of scattered factors up to
length $k \in \N$ is of ongoing research interest. This congruence has attracted recent interest in language
theory \cite{DBLP:journals/tcs/KimHKS23,DBLP:conf/dlt/KimHKS23} and pattern
matching \cite{DBLP:conf/rp/FleischmannKKMNSW23,DBLP:journals/tcs/KimKH24}.
Closely related to Simon's congruence is the $k$-universality of words (also
known as $k$-richness) that receives attention from combinatorial
\cite{DBLP:conf/dlt/BarkerFHMN20,DBLP:conf/stacs/DayFKKMS21,fleischmann2021scattered,DBLP:conf/cwords/SchnoebelenV23}
and language theoretic perspectives \cite{DBLP:conf/isaac/AdamsonFHKMN23,DBLP:conf/isaac/FazekasKMMS24}.
A word is called $k$-universal if it contains all scattered factors of length
$k$. This line of research is continued by characterising words that are not
$k$-universal by their \emph{absent scattered factors}, which are investigated in
\cite{DBLP:journals/corr/abs-2407-18599,DBLP:journals/tcs/FleischmannHHHMN23,DBLP:conf/fct/FleischmannHHN23,DBLP:journals/fuin/KoscheKMS22,DBLP:conf/cocoa/Tronicek23}.

While lossy transmissions and absent scattered factors have been extensively
studied, there is no line of research that allows the transmitted data to be
scrambled \cite{DBLP:journals/tcs/KarhumakiPRW17,DBLP:conf/cwords/Puzynina19}.
A famous application in which the order of letters is not relevant is the
\emph{binary jumbled pattern matching} (BJPM)
\cite{DBLP:journals/ijfcs/BurcsiCFL12,DBLP:journals/algorithmica/GagieHLW15}:  given a word over the binary alphabet $\{0,1\}$
and two natural numbers $x$ and $y$, decide whether the word contains a factor with $x$
zeros and $y$ ones. Here, the order of the occurring letters is irrelevant (hence the
name \emph{jumbled} pattern matching). The description of a word by only the
number of occurrences of each letter is given by its \emph{Parikh vector}
\cite{DBLP:journals/jacm/Parikh66,DBLP:journals/eatcs/Salomaa03}.
Therefore, we introduce and study the new notion of \emph{jumbled} scattered
factors by combining the research on (absent) subsequences and the one of
jumbled (abelian) words. We define a word $u$ to be a \emph{jumbled scattered factor}
of a word $w$ if $u$ can by constructed by the letters of $w$, i.e., $u$ has a
Parikh vector that is componentwise smaller than or equal to the Parikh vector
of $w$. To distinguish those jumbled scattered factors in a more fine-grained
way, we say that a jumbled scattered factor $u$ of a word $w$ has $\ell$ jumbles
if at most $\ell$ letters of $u$ occur at wrong places, i.e., $|u|-\ell$ letters
occur in the correct order as a scattered factor in $w$. For instance, the words
from before $\mathtt{route}$ and $\mathtt{tempo}$ are jumbled scattered factors
of $\mathtt{computer}$ since all their letters occur in $\mathtt{computer}$.
Further, $\mathtt{route}$ is $1$-jumbled since $\mathtt{oute}$ is a scattered
factor of $\mathtt{computer}$; $\mathtt{tempo}$ is $3$-jumbled since
$\mathtt{te}$ or $\mathtt{mp}$ are scattered factors of $\mathtt{computer}$
while the other letters occur in the wrong order. Even though we show that the
minimal number of jumbles of a jumbled scattered factor in a word is closely
related to the longest common scattered factor (lcsf) of those two words, this
work is significantly different than the work on longest common scattered
factors since we are taking a purely combinatorial perspective and partition the
class of jumbled scattered factors by the number of jumbles where lcsf is just one part
of this partition.


{\bf Own Contribution.} In this work, we introduce jumbled scattered factors
and start to explore the broad field of scattered factors in the jumbled case.
First, we have a look at the \emph{jumble index} of a word $u$ in a word $w$ which is the minimal amount of jumbles necessary to make $u$ a scattered
factor of $w$ (\Cref{sec:indexdivision}). 
Then, we continue with a jumbled variant of Simon's
congruence (which is only an equivalence relation) showing that there are finitely many equivalence classes of possible sizes $1,|\Pi(w)|$, or $\infty$ and giving results on the structure of the
classes where $\Pi(w)$ denotes the set of words with the same Parikh vector as $w$ (\Cref{sec:jumblesimon}). In \Cref{sec:universality}, we investigate
jumbled $k$-universality and search for the minimal number of jumbles such that
the jumbled scattered factor set of a given word contains every $k$-length word.


%% file: prelims.tex
Let $\N$ denote the natural numbers starting with $1$ and let $\N_0=\N\cup\{0\}$.
Further, let $[n]=\{1,\ldots,n\}$ for some $n\in\N$ and $[n]_0=[n]\cup\{0\}$.
An {\em alphabet} $\Sigma$ is a finite, non-empty set whose elements are called {\em letters}.
As $\Sigma$ is finite, it can be assumed ordered and w.l.o.g. we take $\Sigma=\{\ta_1,\dots,\ta_{\sigma}\}$
for some $\sigma\in\N$ with $\ta_i<\ta_{i+1}$ for all $i\in[\sigma-1]$.
A finite concatenation of letters from $\Sigma$ is called a {\em word}. The set of words $\Sigma^{\ast}$ with concatenation as binary operation and
with the empty word $\varepsilon$, i.e., the word without any letters, as the neutral element forms a free monoid. For $M\subseteq\Sigma$, set $\cmpl{M}=\Sigma\backslash M$.
The number of letters of a word $w\in\Sigma^{\ast}$ is called its {\em length} and denoted by $|w|$; this implies $|\varepsilon|=0$. Let $\Sigma^k=\{w\in\Sigma^{\ast}|\,|w|=k\}$ and define $\Sigma^{\leq k}$ similarly for all $k\in\N_0$.
The letter of $w\in\Sigma^{\ast}$ at position $i\in[|w|]$ is denoted by $w[i]$. For $w\in\Sigma^{\ast}$ and $\ta\in\Sigma$ define
$|w|_{\ta}=|\{i\in[|w|]|\,w[i]=\ta\}|$, i.e. the number of occurrences of $\ta$ in $w$. Based on this define the {\em Parikh vector} of $w\in\Sigma^{\ast}$ by $p(w)=(|w|_{\ta_1},\ldots,|w|_{\ta_{\sigma}})\in\N_0^{\sigma}$.
We denote the \nth{$i$} entry of $p(w)$ by $p(w)[i]$ for $i \in [\sigma]$. 
It is obvious that two  words with the same Parikh vector have the same length and are permutations of each other like \texttt{flow} and \texttt{wolf}. Given two words $u,v\in\Sigma^{\ast}$, we say $p(u)\leq p(v)$ ($p(u) < p(v)$) if $|u|_{\ta_i}\leq|v|_{\ta_i}$ ($|u|_{\ta_i}<|v|_{\ta_i}$) for all $i\in[\sigma]$. If neither $p(u)\leq p(v)$ nor $p(v)\leq p(u)$ holds, $p(u)$ and $p(v)$ are said to be {\em incomparable}. We extend this notion by saying that $u$ and $v$ are incomparable (w.r.t. their Parikh vector) if $p(u)$ and $p(v)$ are incomparable. We use the maximum norm $\norm{\mathbf{x}}_\infty = \max_{i \in [n]} \abs{x_i}$ for all $\mathbf{x} \in \Z^n$ to refer to the absolute maximum component of a vector.

The main object of interest in our work are {\em scattered factors} also known as subsequences or subwords.

\begin{definition}\label{def:SF}
	Let $w \in \Sigma^*$. A word $u \in \Sigma^*$ is called a \newterm{scattered factor}
	of $w$ if there exists $v_1, ..., v_{\abs{u} + 1} \in \Sigma^*$ such that
	$w = v_1 u[1] ... u[\abs{u}] v_{\abs{u} + 1}$. We denote the set of all scattered factors of $w$ by $\ScatFact(w)$. The largest $k$ such that $\ScatFact(w)\supseteq\Sigma^k$ is called the {\em universality index} of $w$, denoted by $\iota(w)$.
\end{definition}

\begin{definition}\label{crapdef}
	For a set $S\subseteq\Sigma^\ast$ and $k\in\N_0$ define $S_k=S\cap\Sigma^k$ and $S_{\leq k}=\bigcup_{i=0}^kS_i$.
\end{definition}

We mostly use \Cref{crapdef} for $S=\ScatFact(w)$ for a given $w\in\Sigma^{\ast}$.
Roughly spoken, a word $u$ is a scattered factor of a word $w\in\Sigma^{\ast}$, if some of $w$'s letters are deleted and the remaining ones form $u$. For instance, \texttt{allow} is a scattered factor of \texttt{cauliflower} but \texttt{owl} is not since the \texttt{l} occurs only before \texttt{ow}. The latter motivates the definition of {\em jumbled scattered factors}: not the word itself is a scattered factor but a permutation of it is.

\begin{definition}\label{jscatfact}
	Let $w \in \Sigma^*$. A word $u \in \Sigma^*$ is called a \newterm{jumbled scattered factor} of $w$ (denoted $u \jdiv w$) if there exists $v \in \ScatFact(w)$ such that $p(u) = p(v)$. If neither $u\jdiv v$ nor $v\jdiv u$ holds, we say that $u$ and $v$ are {\em incomparable} (w.r.t. jumbling). The set of jumbled scattered factors of $w$ is denoted by $\JScatFact(w)$.
\end{definition}

\begin{remark}\label{incomparable}
	Two words are incomparable w.r.t.\ their Parikh vector if and only if they are incomparable w.r.t.\ jumbling. Hence, in the following we will not distinguish between both incomparabilities.
\end{remark}

Thus, in \texttt{cauliflower}, not only \texttt{owl} but also \texttt{wolf}, \texttt{recall}, \texttt{failure}, and \texttt{firewall} are jumbled scattered factors. Notice that \Cref{jscatfact} implies that every scattered factor is also a jumbled scattered factor and $\JScatFact(w)$ is finite for all $w\in\Sigma^{\ast}$. Inspecting \texttt{recall} and \texttt{wolf} in more detail one notices that the longest scattered factor which \texttt{recall} and \texttt{cauliflower} have in common (\texttt{call}) is longer than the one \texttt{wolf} and \texttt{cauliflower} have in common (\texttt{lf}). Before we relate the well known concept of the {\em longest common scattered factor} and the jumbled scattered factors, we introduce some definitions for the number of required jumbles.

\begin{definition}\label{def:jidx}
	Let $w \in \Sigma^*$. The word $u \in \Sigma^*$ is called a \newterm{scattered factor
		of $w$ with $\ell$ jumbles} (or $u$ is $\ell$-jumbled in $w$) for $\ell \in [\abs{u}-1]_0$ if $p(u)\leq p(w)$ and there exists $v\in\ScatFact_{\abs{u} - \ell}(u)$ such that $v \in \ScatFact(w)$ holds.
	The set of scattered factors of $w$ with $\ell$-jumbles is denoted by $\JScatFact(w,\ell)$.
\end{definition}


Note that, for $w\in\Sigma^\ast$, $u \in \JSF(w, \ell)$ implies that every subsequence of $u$ is at most $\ell$-jumbled in $w$. By \Cref{def:jidx} we have that \texttt{owl} is $1$-jumbled and \texttt{wolf} is not only $2$-jumbled but also $3$-jumbled depending on which \texttt{l} in \texttt{cauliflower} is chosen. Having a measure in mind that a jumbled scattered factor is {\em better}, i.e., closer to being a scattered factor, the less jumbles are required (for example in DNA-sequencing or in data transmission), we fix the smallest possible number of needed jumbles.

\begin{definition}\label{jumbleindex}
	Let $w \in \Sigma^*$ and $u \in \JScatFact(w)$. We define the
	\newterm{jumble index of $u$ w.r.t. $w$} by $\delta_w(u) = \min\{\ell \in \N_0 \mid u\in\JScatFact(w,\ell)\}$.
	Let $\SJSF(w,\ell)=\{u\in\JSF(w,\ell)\mid\delta_w(u)=\ell\}$.
\end{definition}

\Cref{jumbleindex} gives us $\delta_{\mathtt{cauliflower}}(\mathtt{wolf})=2$ and $\delta_{\mathtt{cauliflower}}(\mathtt{failure})=3$. Note that $w \in \JScatFact(w, \ell)$ for all $\ell \in [\abs{w}-1]_0$ but $w \in \SJSF(w, \ell) \iff \ell = 0$ for $w \in \Sigma^*$.

\begin{remark}\label{remark:0jumbles}
	For $w \in \Sigma^*, k \in \N$, every $u \in \JumbledScatFact_k(w)$ with
	$\delta_w(u) > 0$ is of length at least $\iota(w) +1$ since $|u|\leq\iota(w)$ implies $u\in\JScatFact(w,0)=\ScatFact(w)$.
	The same argument gives us that unary jumbled scattered factors $u$ yield $\delta_w(u)=0$.
\end{remark}

\begin{definition}
	Let $u, v \in \Sigma^*$. We define the set of {\em common scattered factors} by
	$\CSF(u, v) = \ScatFact(u) \cap  \ScatFact(v)$ and $\LCSF(u,v)=\{x\in\CSF(u,v)\mid \forall y\in\CSF(u,v):\,|x|\geq |y|\}$ as the set of {\em longest common scattered factors of $u$ and $v$}. Their length is denoted by $\lcsf(u, v)$.
\end{definition}

Some relations between the $\LCSF$ and $\JScatFact$ are immediate.

\input{lcs}

Based on Simon's congruence \cite{DBLP:conf/automata/Simon75}, we also define a variant for jumbled scattered factors.

\begin{definition}
	Let $v,w \in \Sigma^*$ and $k,\ell \in \N$. We say that $v$ is \emph{$\ell$-jumble Simon $k$-equivalent} to $w$, for short $v \jsimkl w$ iff $\JScatFact_{\leq k}(v, \ell) = \JScatFact_{\leq k}(w, \ell)$. We denote the equivalence class of $w$ under $\jsimkl$ by $[w]_{k,l} = \{v\in\Sigma^\ast\mid v\jsimkl w\}$.
\end{definition}

\begin{lemma}
	For $k,\ell \in \N, w \in \Sigma^*$, $\JScatFact_{k}(u, \ell) = \JScatFact_{k}(v, \ell)$ implies $\JScatFact_{\leq k}(u, \ell) = \JScatFact_{\leq k}(v, \ell)$.
\end{lemma}
\begin{proof}
	Suppose $\JScatFact_{k}(u, \ell) = \JScatFact_{k}(v, \ell)$ and $\JScatFact_{\leq k}(u, \ell) \neq \JScatFact_{\leq k}(v, \ell)$. Thus, w.l.o.g., there exists a $k' < k$ (choose $k'$ as the largest) and $x \in \Sigma^*$ with $x \in \JScatFact_{k'}(u, \ell)$ (in particular $\abs x = k'$) and $x \notin \JScatFact_{k'}(v, \ell)$. Now there exist $y_1, \ldots, y_{k'+1} \in \Sigma^*$ of length $|y_1 \cdots y_{k'+1}| = k-k'$ such that $y_1 x[1] y_{2} x[2] y_{3} x[3] \cdots y_{k'} x[k'] y_{k'+1}\in \JScatFact_k(u,\ell)$. Thus, $y_1 x[1] y_{2} x[2] y_{3} x[3] \cdots y_{k'} x[k'] y_{k'+1} \in \JScatFact_k(v,\ell)$ which implies a contradiction: $x \in \JScatFact_{k'}(v,\ell)$.
	\qed\end{proof}

This leads to the problem of investigating the jumbled variant $\jsimkl$ of Simon's congruence relation.

\begin{problem}\label{problem:simon}
	Given two words $v,w \in \Sigma^*$ and two integers $k,\ell$ is $v \jsimkl w$?
\end{problem}

Tightly related to the classical Simon's congruence with scattered factors is the notion of universality \cite{DBLP:conf/dlt/BarkerFHMN20,fleischmann2021scattered,DBLP:conf/cwords/SchnoebelenV23} that we extend to the jumbled scattered factors.
Intuitively, one may think that the notion of $k$-universality can be extended to $\ell$-jumble $k$-universality by requiring that a word contains all words of length $k$ as scattered factors with exactly $\ell$-jumbles, i.e., $\SJSF_k(w,\ell) = \Sigma^k$. This definition would invoke that no word is $\ell$-jumble $k$-universal for $\ell>0$: all scattered factors of length $k$ are scattered factors with exactly $0$ jumbles and thus are missing in $\SJSF_k(w,\ell)$ (cf. \Cref{remark:0jumbles}). For this reason we relax the generalisation of universality to jumbled-universality by allowing {\em up to} $\ell$ jumbles. This idea is captured in the following definition.

\begin{definition}
	A word $w \in \Sigma^*$ is called \emph{$\ell$-jumble $k$-universal} if $$\bigcup_{\ell' \in [0, \ell]} \SJSF_k(w,\ell') = \Sigma^k.$$
\end{definition}

This leads to the problem of determining the smallest value for $\ell$  to obtain universality by the union of all $\ell$-jumbled scattered factors.

\begin{problem}\label{problem:univ}
	Given $w\in\Sigma^{\ast}$ and $k\in\N$ with $\iota(w)<k$. Determine the minimal $\ell\in\N$ such that $\bigcup_{\ell'\in[\ell]_0}\SJSF_k(w,\ell')=\Sigma^k$.
\end{problem}

We conclude with a tool that has proven useful for the study of scattered factor universality: the arch factorisation that was defined by Hébrard \cite{DBLP:journals/tcs/Hebrard91}.

\begin{definition}
	For a word $w \in \Sigma^*$ the \emph{arch factorisation} is given by
	$w = \ar_1(w)\cdots\ar_k(w)r(w)$ for $k \in \N_0$ with\\
		- $\iota(\ar_i(w)) = 1$ for all $i \in [k]$,\\
		- $\ar_i(w)[\abs{\ar_i(w)}] \notin \letters(\ar_i(w)[1..\abs{\ar_i(w)} - 1])$, and\\
		- $\letters(r(w)) \subset \Sigma$.\\
	The unique last letters of each arch $m(w) = \ar_1(w)[\abs{\ar_1(w)}]\cdots\ar_k(w)[\abs{\ar_k(w)}]$ build the
	\emph{model} of $w$.
\end{definition}

For a given word $w\in\Sigma^{\ast}$ the arch factorisation decomposes a word into uniquely defined arches: read $w$ from left to right and every time the complete alphabet was read, close an arch. For instance, consider $w=\mathtt{alfalfa}$
we have the decomposition $\mathtt{alf}.\mathtt{alf}.\ta$ where the dots indicate the end of an arch.

%% file: lcs.tex

\ifpaper
\input{lcsfparakh_to_jsf}
\else
\input{lcsfparakh_to_jsf}
\input{PROOFlcsfparakh_to_jsf}
\fi

%
%
%
%
%
%

%% file: lcsfparakh_to_jsf.tex
\begin{lemma}\label{lcsfparakh_to_jsf}
(1) Let $u, w \in \Sigma^*$, $\ell \in [\abs{u}-1]_0$ and $v \in \CSF(u, w)$. If
	$\ell + \abs{v} = \abs{u}$ and $p(u) \leq p(w)$, then $u\in\JSF(w,\ell)$.\\
(2) For $w\in\Sigma^{\ast}$ and $u\in\JScatFact(w,\ell)$  there exists $v \in \CSF(u, w)$ such that
	$l + \abs{v} = \abs{u}.$\\
(3) 	For $w \in \Sigma^*$, $\ell\in\N$, we have $u\in \SJSF(w,\ell)$ iff $\ell = \abs{u} - \lcsf(u, w)$. Thus, $\delta_w(u) = \abs{u} - \lcsf(u, w)$.
\end{lemma}

%% file: PROOFlcsfparakh_to_jsf.tex
\begin{proof}
	(1) follows immediately by $v \in \ScatFact(u)$. 	By definition of $u$, there exists $v \in \ScatFact(w)$ with
	$\abs{v} = \abs{u} - l$ and $v\in\ScatFact(u)$. Thus,
	$v \in \CSF(u, w)$ and $l + \abs{v} = \abs{u}$ holds which proves (2). If we choose a longest common scattered factor, we get immediately the duality between
the jumble index and longest common scattered factors. This finishes the proof of (3).\qed
\end{proof}

%% file: indexdivision.tex
We start this section with some basic properties of the relation $\jdiv$ (in this section we write $u \jdiv w$ instead of $u \in \JScatFact(w)$ since we are studying the relation). It is immediate that $|_J$ is reflexive and transitive. In contrast to the scattered factor relation \cite{lothaire}, $\jdiv$ is neither symmetric nor antisymmetric: from $u\jdiv v$ and $v\jdiv u$ it only follows $p(u)=p(v)$. Let $u,u',w,w' \in \Sigma^*$. If $u \jdiv w$ and $u' \jdiv w'$ then $uu' \jdiv ww'$. The relation $\jdiv$ is a preorder on $\Sigma^*$ that is compatible with $\cdot$ on $\Sigma^*$.


\Cref{incomparable} motivates the following proposition, a similar results was shown for scattered factors in \cite[Theorem 6.1.2]{lothaire}. 
The proof follows immediately by \cite{KRUSKAL1972297}.

\input{pairwiseincomparable}


A word $u$ can be a scattered factor of $w$ with a different number of jumbles.
For example, consider $w = \mathtt{abba}$ and $u = \mathtt{ba}$. Here, the number of jumbles can
be either zero or one.
The discrete interval $[\delta_w(u), \abs{u}-1]$ describes the set of all possible
numbers of jumbles for which $u$ is a jumbled scattered factor of $w$. By
\Cref{lcsfparakh_to_jsf} we can choose any common scattered factor of $u$ and
$w$ as long as we choose the number of jumbles appropriately. Or to express it
differently, once we have determined $\delta_w(u)$, we can simply remove a
character from a fixed longest common scattered of both words, thus increasing
the number of jumbles that are used for $u$ to be a jumbled scattered factor.

\ifpaper
\input{jscatfactspectrum}

\input{PROOFjscatfactspectrum}
\else

\input{jscatfactspectrum}
\input{PROOFjscatfactspectrum}
\fi

This statement also holds for the length-restricted jumbled scattered factors.

\ifpaper
\input{lemjsfinclusion}

\input{PROOFlemjsfinclusion}
\else

\input{lemjsfinclusion}
\input{PROOFlemjsfinclusion}
\fi

Whereas in the field of scattered factors, we have that $u\in\ScatFact(w)$ and $w\in\ScatFact(u)$ implies $w=u$, in the field of jumbled scattered factors, it only implies that $u$ and $w$ are permutations of each other. This leads to the following result regarding the number of jumbles.

\ifpaper
\input{equaldelta}

\input{PROOFequaldelta}
\else

\input{equaldelta}
\input{PROOFequaldelta}
\fi

Another difference to {\em unjumbled} scattered factors is the concatenation of jumbled scattered factors.
The words $\mathtt{life}$ and $\mathtt{low}$ are both scattered factors of $w=\mathtt{cauliflower}$ and even the sum of their Parikh vectors is smaller than the Parikh vector of $w$ but $\mathtt{lowlife}$ is not a scattered factor of $w$. In the scenario of jumbled scattered factors,
we have that $\delta_{\mathtt{cauliflower}}(\mathtt{low\cdot life})=3$, thus the jumble index is higher than the sum of the jumble indices. The sum is the lower bound, e.g. $\delta_{\mathtt{cauliflower}}(\mathtt{ufo\cdot war})=1$.

\ifpaper
\input{uvinw}

\input{PROOFuvinw}
\else

\input{uvinw}
\input{PROOFuvinw}
\fi

Moreover, for the concatenation of not only the jumbled scattered factors but also the words,
we can bind the jumble index. Consider for instance $w = \mathtt{cba}$, $u=\mathtt{abc}
\in\JScatFact(w)$, with $\delta_w(u)=2$, and $w' = \mathtt{bcbaa}$, $v=\ta\tb\in\JScatFact(w')$, with $\delta_{w'}(v)=1$. In $ww'=\mathtt{cbabcbaa}$, we have $u\ta$ as a scattered factor and thus $\delta_{ww'}(uv) = 1$.
For the total number of jumbles, we had the advantage that $u$ is a scattered factor of $ww'$.
Generally, in this case we have the sum $\delta_w(u)+\delta_{w'}(v)$ as upper bound for $\delta_{ww'}(uv)$.

\ifpaper
\input{uvinww}

\input{PROOFuvinww}
\else

\input{uvinww}
\input{PROOFuvinww}
\fi

We finish this section an immediate result on the relation of the index of a word and its index when extending the word by one letter.

\begin{remark}\label{lemma:append}
	Let $w \in \Sigma^*$ and $u \in \JSF(w)$.
	Then, $\delta_w(ux) - \delta_w(u) \leq 1$ holds for all $x \in \Sigma$ with $p(ux) \leq p(w)$.
\end{remark}

%

%% file: pairwiseincomparable.tex
\begin{proposition}\label{pairwiseincomparable}
	Let $S \subseteq \Sigma^*$ be a set of words that are not pairwise comparable with $\jdiv$. Then $S$ is finite. 
\end{proposition}

%% file: jscatfactspectrum.tex
\begin{corollary}\label{jscatfactspectrum}
  Let $w \in \Sigma^*$ and $u \in \JScatFact(w)$. Then,
  $u \in \JScatFact(w, \ell)$ for all $\ell \in [\delta_w(u), \abs{u}-1]$.
\end{corollary}

%% file: PROOFjscatfactspectrum.tex
\begin{proof}
Let $u\in\JScatFact(w)$, thus we have  $p(u)\leq p(w)$. Then there exists $v\in\ScatFact(w)$ with $p(u)=p(v)$. All $v'\in\CSF(u,v)$ provide by \Cref{lcsfparakh_to_jsf} an $\ell$ with $u\in\JScatFact(w,\ell)$. The claim follows by $0\leq|v'|\leq\lcsf(u,v)$.\qed
\end{proof}

%% file: lemjsfinclusion.tex
\begin{lemma}\label{lem:jsf_inclusion}
  Let $w \in \Sigma^*$, $k\leq \abs{w}$, and $\ell\in[k-2]$. Then, $\JScatFact_k(w, \ell) \subseteq \JScatFact_k(w, \ell+1)$ and $\JScatFact_k(w, i) = \emptyset$ for all $i \geq k$.
\end{lemma}

%% file: PROOFlemjsfinclusion.tex
\begin{proof}
  Let $\ell \in [0, k - 2]$ and $u \in \JScatFact_k(w, \ell)$. There
  exists $v \in \CSF(u, w)$ such that $\ell + \abs{v} = \abs{u}$. Then, $v\neq\varepsilon$
  since $\abs{v} = \abs{u} - \ell \geq 1$ holds. We can remove
  one character from $v$ and get $v' = v[1..\abs{v} - 1] \in \CSF(u, w)$.
  Thus, $u$ is $(\ell + 1)$-jumbled w.r.t. $w$, since
  $\ell + 1 + \abs{v'} = \abs{u}$ holds. Furthermore, since all quantities must
  be non-negative in $\ell + \abs{v} = \abs{u} = k$, we cannot choose
  $\abs{v} < 0$ to accommodate for $\ell > k$.
    Suppose that $\JScatFact(w, \ell) \neq \emptyset$ for some $\ell > \abs{w}$,
  then there exists $u \in \JScatFact(w, \ell)$ and thus a $v \in \CSF(u, w)$
  with $\abs{v} + \ell = \abs{u}$. However, $\abs{w} < \ell = \abs{u} - \abs{v}$
  contradicts $\abs{u} \leq \abs{w}$. This concludes the proof.\qed
\end{proof}

%% file: equaldelta.tex
\begin{lemma}\label{equaldeltas}
  Let $u,w \in \Sigma^*$ such that $u\jdiv w$ and $w\jdiv u$. Then,
  $|u|=|w|$,  $\delta_w(u) = \delta_u(w)$.
\end{lemma}

%% file: PROOFequaldelta.tex
\begin{proof}
The first claim follows immediately by the definition, in particular $p(u)\leq p(w)$ and
$p(w)\leq p(u)$ imply $p(u)=p(w)$ which implies $\abs u = \abs w$. The second part follows from $\lcsf(u,w) = \lcsf(w,u)$: $\delta_u(w) = \abs w-\lcsf(w,u) = \abs u-\lcsf(u,w) =\delta_w(u)$. \qed
\end{proof}

%% file: uvinw.tex
\begin{lemma}\label{uvinw}
	Let $w \in \Sigma^*$ and $uv \in \JumbledScatFact(w)$. Then,
	$\delta_w(uv) \geq \delta_w(u) + \delta_w(v)$.
\end{lemma}

%% file: PROOFuvinw.tex
\begin{proof}
Choose $\ell\in\N_0$ minimal with $uv\in\JScatFact(w,\ell)$. Then, we have $p(uv)\leq p(w)$ and there exists 
\[
x\in\ScatFact_{|uv|-\delta_w(uv)}(uv)\cap\ScatFact(w).
\]
Thus, there exist $x_u,x_v\in\Sigma^{\ast}$ with $x=x_ux_v$, $x_u\in\ScatFact(u)$, and $x_v\in\ScatFact(v)$. Moreover, we have
\[
|x_u|+|x_v|=|x_ux_v|=|x|=|uv|-\delta_w(uv)=|u|+|v|-\delta_w(uv).
\]
Thus, there exist $\ell_u$ and $\ell_v$ with $|x_u|=|u|-\ell_u$, $|x_v|=|v|-\ell_v$ and $\delta_w(uv)=\ell_u+\ell_v$. Accordingly, $x_u\in\ScatFact_{\abs u - \ell_u}(u)$ and $x_v\in\ScatFact_{\abs v - \ell_v}(v)$. Because $x$ is a subsequence of $w$, $x_u$ and $x_v$ also are subsequences of $w$. It follows $\delta_w(u)\leq\ell_u$, $\delta_w(v)\leq\ell_v$ and, thereby, $\delta_w(uv) = \ell_u+\ell_v\geq\delta_w(u)+\delta_w(v)$.\qed
\end{proof}

%% file: uvinww.tex
\begin{lemma}\label{uvinww}
	Let $u,v,w,w' \in \Sigma^*$. If $u \in \JSF(w)$ and $v \in \JSF(w')$ then $uv \in \JSF(ww')$  with $\min\{\delta_{w}(u),\delta_{w'}(v)\} \leq \delta_{ww'}(uv) \leq \delta_{w}(u) + \delta_{w'}(v)$.
\end{lemma}

%% file: PROOFuvinww.tex
\begin{proof}
	The fact that $uv$ is a jumbled scattered factor of $ww'$ follows directly from Definition~\ref{jscatfact}.
	The inequality holds by the minimality of the jumbled index.
	\qed
\end{proof}

%% file: jumblesimon.tex
In this section, we investigate the extension of the famous Simon congruence \cite{DBLP:conf/automata/Simon75} in the jumbled scenario.
Note that there are subsets $M\subset \Sigma^{\leq k}$ such that for all $\ell\in\N_0$ and $S\subseteq\Sigma^\ast$ it holds $M\neq\bigcup_{w \in S} \JScatFact_{\leq k}(w,\ell)$ (this property is already true for $\ell =0$ , i.e., scattered factors). For instance, there exists no $w \in \Sigma^*, \ell\in\N_0$ such that
$\JScatFact_{\leq k}(w,\ell)=\{\ta, \tb\tb\}\subset \Sigma^{\leq 2}$ -
at least one of $\ta \tb$ or $\tb \ta$ needs to be a ($0$-jumbled) scattered factor. Moreover,
if $|w| < k$ then $[w]_{k,\ell} = \{w\}$. Thus, we consider only values $\iota(w) \leq k \leq |w|$ and $0 \leq \ell < k$.
We extend results from \cite{lothaire,DBLP:conf/automata/Simon75} on $\sim_k$ to the jumbled case $\jsimkl$ and continue with new results for the structure of jumbled Simon's equivalence classes (Problem~\ref{problem:simon}).

Before we start with insights about $\jsimkl$, we present some examples for the kinds of equivalence classes that may occur. First, consider some singleton classes $[\mathtt{ab}]_{\jsim_{2,0}}$ or $[\mathtt{abcd}]_{\jsim_{3,2}}$. Second, the class $[\mathtt{ab}]_{\jsim_{2,1}} = \{\mathtt{ab}, \mathtt{ba}\}$ containing exactly the two words with the Parikh vector $(1,1)$. And the infinite class $[\mathtt{aab}]_{\jsim_{2,1}} = \{\ta^n\tb\ta^m\mid n,m \in \N_0,n+m\geq 2\}$ that all share the $1$-jumbled scattered factor set of length $2$, i.e., $\JSF_2(w,1) = \{\mathtt{aa}, \mathtt{ab}, \mathtt{ba}\}$ for all $w\in [\mathtt{aab}]_{\jsim_{2,1}}$.


\ifpaper
\input{congruencerelation}
\input{PROOFcongruencerelation}

\else

\input{congruencerelation}
\input{PROOFcongruencerelation}
\fi

Note that $\jsimkl$ is in contrast to $\sim_k$ not a congruence relation. As an example consider $w = \mathtt{cd}$ and $v = \mathtt{dc}$. We have $w \jsim_{3,1} v$ but $\mathtt{ab}w \not \jsim_{3,1} \mathtt{ab}v$ since $\mathtt{cda} \in \JScatFact_3(\mathtt{ab}w,1)$ and $\mathtt{cda} \notin \JScatFact_3(\mathtt{ab}v,1)$.

\ifpaper
\input{jsimbasics}
\input{PROOFjsimbasics}
\else

\input{jsimbasics}
\input{PROOFjsimbasics}
\fi

\ifpaper
\input{jsimkl_letter_ext}
\input{PROOFjsimkl_letter_ext}
\else

\input{jsimkl_letter_ext}
\input{PROOFjsimkl_letter_ext}
\fi

We now continue by a characterisation of the equivalence classes of $\jsimkl$.
Note that $\ell = 0$ implies by \cite{lothaire} that $\abs{[w]_{\jsim_{k,\ell}}} \in \{1,\infty\}$ since  $\jsim_{k,0}$ is the same as Simon's congruence relation. Nevertheless, the classes of infinite size differ if $\ell \geq 1$, e.g., $[\mathtt{aab}]_{\sim_3}$ is a singleton congruence class whereas $[\mathtt{aab}]_{\jsim_{2,\ell}}$ is an infinite equivalence class for $\ell \geq 1$.
One may think that singleton classes of $\jsimkl$ only exist for $\ell=0$, i.e.,   in the \emph{classical} Simon's congruence. 
Consider $w=\mathtt{abcd}\in\{\ta,\tb,\tc,\td\}^{\ast}$ providing a counterexample: we
have $|\JScatFact_3(\mathtt{abcd},1)|=20$ (missing the words
$\tc\tb\ta, \td\tb\ta, \td\tc\ta, \td\tc\tb$); no word of length $>4$ is in this class because $\tx\tx\ty\notin\JScatFact_3(w,1)$, for all $\tx,\ty\in\Sigma$, and similarly any
word missing one of the four letters cannot be equivalent since we would lose jumbled scattered factors; lastly no permutation
of $\ta\tb\tc\td$ is in the class since there exists a jumbled scattered factor of any permutation of $\mathtt{abcd}$ that is contained in the set $\JScatFact_3(\mathtt{abcd},1)$ but needs more than $1$ jumble.

	%

\begin{remark}
	Let $k,\ell \in \N$. Note, that two words $w,w' \in \Sigma^*$ with $\al(w) = \al(w')$ and $|w|_\ta <k$ for all $\ta \in \Sigma$ cannot be equivalent if there exists $\tb \in \Sigma$ such that $|w|_\tb \neq |w'|_\tb$.
\end{remark}

As a first step we consider infinite classes. The main idea on the characterisation of infinite equivalence classes is that a letter that already occurs more than or equal to $k$ times in a sufficiently large block can also occur more often: all further occurrences are superfluous for $k$-length equivalence. Recall that we only need to investigate $0\leq \ell<k$.

\ifpaper
\input{char_of_infinite_classes}
\input{PROOFchar_of_infinite_classes}
\else

\input{char_of_infinite_classes}
\input{PROOFchar_of_infinite_classes}
\fi

Note that the restriction on the length of the $\ta$-block is necessary. Consider for example $k = 6, \ell = 1$, $w = \mathtt{cb}\ta^3\mathtt{cb}\ta^3$ and $w' = \mathtt{cb}\ta^4\mathtt{cb}\ta^3$. One can verify that $\ta^4\mathtt{bc} \in \JScatFact_6(w,2)$ with $\delta_w(\ta^4\mathtt{bc}) = 2$ and $\ta^4 \mathtt{bc} \in \JScatFact_6(w',1)$. 
Now, we show that if there exist two equivalent words with a different Parikh vector then there always exists a word within the same equivalence class that contains the two equivalent words as its jumbled scattered factors. Following the idea given in \Cref{prop:char_of_infinite_classes} this word will be constructed by the addition of superfluous letters that occur more than $k$ times anyway.

\ifpaper
\input{superjumblesf}
\input{PROOFsuperjumblesf}
\else

\input{superjumblesf}
\input{PROOFsuperjumblesf}
\fi

Note that we can also conclude $w,u \in \JScatFact(v)$ from the above corollary if $p(w) \leq p(v)$ or $p(u) \leq p(v)$ additionally holds for $v$. 

Now we show when two equivalent words are abelian equivalent.

\ifpaper
\input{parikhequal}
\else
\input{parikhequal}
\input{PROOFparikhequal}
\fi

\begin{corollary}
	Let $w \in \Sigma^*, k \in \N_{>1}, \ell \in \N$ with $|w|_\ta < k$ for all $\ta \in \Sigma$. We have $[w]_{\jsim_{k,\ell}} = \Pi(w)$ iff $\ell = k-1$.
\end{corollary}

	%
	%
	%

%

Whether other sizes of equivalence classes are possible stays an open problem for future research.

%% file: congruencerelation.tex
\begin{lemma}\label{congruencerelation}
	$\jsimkl$ is a equivalence relation on $\Sigma^*$ of finite index.
\end{lemma}

%% file: PROOFcongruencerelation.tex
\begin{proof}
It is immediate that $\jsimkl$ is an equivalence relation.

	
	Further, the index of $\jsimkl$ is equal to the number of distinct suitable choices for $\JScatFact_{\leq k}(w,\ell)$ for all $w \in \Sigma^*$ (each equivalence class is uniquely determined by the set of its jumbled scattered factors). 
	We have 
  
	\[
    \left|\bigcup_{w \in \Sigma^*} \JScatFact_{\leq k}(w,\ell)\right| \leq |\mathcal P (\Sigma^{\leq k})| = 2^{|\Sigma^{\leq k}|}.\quad\qed
	\]
\end{proof}

%% file: jsimbasics.tex
\begin{lemma}\label{lem:jsimbasics}
	Let $\ell,k\in\mathbb N$, $\ell\leq k-2$ and $w,v \in \Sigma^*$ satisfy $w \jsim_{k+1,\ell} v$. Then $w \jsim_{k+1,\ell+1} v$, $w \jsim_{k,\ell} v$, and $w \jsim_{k,\ell+1} v$ also hold.
\end{lemma}

%% file: PROOFjsimbasics.tex
\begin{proof}
  First, let $w \jsim_{k+1,\ell} v$ and assume $w \not\jsim_{k+1,\ell+1} v$. Then, without loss of generality, there exists $u\in\Sigma^{m}$, with $m\leq k+1$, such that $\delta_v(u)= \ell+1<m$ and $\delta_w(u)\geq \ell+2$. Let $u=u_1\cdots u_m$, $\hat u=u_{i_1}\cdots u_{i_{m-\delta_v(u)}}$ be a (longest) common scattered factor of $u$ and $v$ and $j\notin\{i_1,\ldots,i_{m-\delta_v(u)}\}$. Then $\tilde u=u_1\cdots u_{j-1}u_{j+1}\cdots u_m$ and $v$ have the same common scattered factor $\hat u$. Accordingly, $\tilde u\in\JScatFact_{\leq k+1}(v,\ell)=\JScatFact_{\leq k+1}(w,\ell)$ and there exists a common scattered factor $\bar u$ of length $\abs{\tilde u} - \ell$ of $\tilde u$ and $w$. Lastly, $\bar u$ is also a common scattered factor of $u$ and $w$ which implies $\delta_w(u)\leq \abs u -\abs{\bar u}=\abs{\tilde u}+1-\abs{\bar u} = \ell+1$.

	Second, $w \jsim_{k,\ell} v$ holds because any jumbled scattered factor of length $k$ of $v$ and, respectively, $w$ occurs as scattered factor in a jumbled scattered factor of length $k+1$ of $v$ and, respectively, $w$. 
	
	The last statement follows as a combination of the previous two.

	Let $v \in \JScatFact(w)$. First, assume that $v \jsimkl w$ holds and suppose that there exists a $u \in \Sigma^*$ with $|u| \leq k$ and $\ell \geq \delta_w(u)$ such that $u \notin \JScatFact(v)$. By definition this implies $v \notin \JScatFact(w)$, a contradiction.
	Second, assume that for all $u \in \Sigma^*$, $|u| \leq k$ and $\ell \geq \delta_w(u)$ implies $u \in \JScatFact(w)$. We are going to show that $v \jsimkl w$ by supposing the contrary. Suppose $v \not \jsimkl w$. Since $v \in \JScatFact(w)$, there exists  $u \in \JScatFact_{\leq k}(w)$  and $\ell \geq \delta_w(u)$ such that $u  \notin \JScatFact(v)$, again a contradiction to the assumption.\qed
\end{proof}

%% file: jsimkl_letter_ext.tex
\begin{lemma}\label{lem:jsimkl_letter_ext}
	For all $u \in \Sigma^*$ and $\ta \in \Sigma$ there exist $p, \ell \in \N_0$, $s \in \Sigma^*$ such that $u \jsim_{p,\ell} u \ta,  s \in \JScatFact_p(u)$ with $\delta_u(s) \leq \ell$ and $s\ta \notin \JScatFact(u)$.  
\end{lemma}

%% file: PROOFjsimkl_letter_ext.tex
\begin{proof}
	The existence of such $p$ follows from \lemref{lem:jsimbasics} since $u \jsim_{0,\ell} u\ta$ and $u \not\jsim_{|u|+1,\ell} u\ta$ hold for every $\ell \in \N_0$. Thus, we may choose the greatest $p$ with $u \jsim_{p,\ell} u\ta$ and $u \not\jsim_{p+1,\ell} u\ta$. Then the existence of $s \in \JScatFact_p(u)$ with $s\ta \notin \JScatFact_p(u)$ follows.\qed
\end{proof}

%% file: char_of_infinite_classes.tex
\begin{proposition}\label{prop:char_of_infinite_classes}
  Let $w \in \Sigma^*$ and $k,\ell \in \N$.
	If there exist $\ta\in\al(w)$, $k'\geq k-\ell$, and $x,y\in\Sigma^{\ast}$ with $|w|_{\ta}\geq k$ and $w=x\ta^{k'}y$
	then $|[w]_{\jsim_{k,\ell}}| = \infty$.
\end{proposition}

%% file: PROOFchar_of_infinite_classes.tex
\begin{proof}
	Let $\ta \in \al(w), k' \geq k-\ell$ and $x,y \in \Sigma^*$ with $|w|_\ta \geq k$ and 
	 $w = x \ta^{k'} y$.
	We now want to show that $w \jsimkl w'$ for $w' = x \ta^{k'+1} y$, i.e., $\JScatFact_k(w,\ell) = \JScatFact_k(w',\ell)$.
	
	Let $u \in \JScatFact_k(w,\ell)$. Thus, there exists $v \in \ScatFact_{k-\ell}(u)$ with $v \in \ScatFact(w)$.  By construction, we have that $v \in \ScatFact(w')$. By $p(u)\leq p(w)\leq p(w')$, we have  $u \in \JScatFact_k(w',\ell)$.
	
	For the other inclusion, let $u' \in \JScatFact_k(w',\ell)$. So, there exists $v' \in \ScatFact_{k-\ell}(u')$ with $v' \in \ScatFact(w')$. Since $k' +1 > k-\ell$ we obtain $p(u')\leq p(w)$, i.e., that at least one of the $\ta$'s in the block of length $k'$ must not be used in $v'$. Thus, we have $v' \in \ScatFact_{k-\ell}(w)$. 
	
	\medskip
	
	The claim follows by induction.\qed
\end{proof}

%% file: superjumblesf.tex
\begin{corollary}\label{prop:superjumblesf}
	For all $w,u \in \Sigma^*,k,\ell\in\N$ with $p(w) \neq p(u)$ and all $\ta_i \in \Sigma$ with $p(w)[i],$ $p(u)[i] \geq k$, $\ta_i^{k'} \in \Fact(w) \cap \Fact(u)$ for $k' \geq k-\ell$, and $w \jsimkl u$ there exists $v \in \Sigma^*\setminus\{u,w\}$ such that  $v \jsimkl w, v \jsimkl u$.
\end{corollary}

%% file: PROOFsuperjumblesf.tex
\begin{proof}
Let $w \jsimkl u$ with $p(w) \neq p(u)$. 
First, we want to show that for all $i \in [\sigma]$ with $p(w)[i] \neq p(u)[i]$ it follows that $p(u)[i],p(w)[i] \geq k$. For the sake of contradiction suppose the contrary, i.e., that there exists $i \in [\sigma]$ with $p(w)[i] \neq p(u)[i]$ and with $p(w)[i] < k$ or $p(u)[i] < k$.
W.l.o.g. let $p(w)[i] < p(u)[i]$. Thus, there exists $x \in \JSF_{k}(u,\ell)$ with $p(x)[i] = p(u)[i]$. Since $p(w)[i] < p(u)[i]$ and at least one of those values is smaller than $k$ we obtain $x \notin \JSF_{k}(w,\ell)$ - a contradiction to $w \jsimkl u$.

Thus, for all $i \in [\sigma]$ with $p(w)[i] \neq p(u)[i]$ we have that $p(w)[i] \geq k$ and $p(u)[i] \geq k$.
Let $\ta_i$ be a letter for which the Parikh vectors of $u$ and $w$ are different. W.l.o.g. let $w$ be the word that contains a larger block of $\ta_i$. By assumption we can write $w$ by $w = x \ta_i^{k'} y$ for some $x,y \in \Sigma^*$ and $k' \geq k-\ell$.
Now we can construct a $v$ following the proof of \Cref{prop:char_of_infinite_classes}:
Let $v = x \ta_i^{k'+1} y$. Thus, $v \jsimkl w \jsimkl u$. It is left to show that $w,u \in \JScatFact(v)$. The fact that $w \in \JScatFact(v)$ is immediately visible. 
		\qed
\end{proof}

%% file: parikhequal.tex
\begin{proposition}\label{parikhequal}
	Let $w, w' \in \Sigma^*, k \in \N_{>1}$ with $|w|_\ta,|w'|_\ta < k$ for all $\ta \in \Sigma$. We have $w \jsim_{k,k-1} w'$ iff $p(w) = p(w')$.
\end{proposition}

%% file: PROOFparikhequal.tex
\begin{proof}
	We have 
	\begin{align*}
		&w  \jsim_{k,k-1} w'\\
		\iff& \{u \in \Sigma^k \mid (\exists v \in \ScatFact_1(u) : v \in \ScatFact(w)) \land p(u)\leq p(w)\} \\
		&= \{u' \in \Sigma^k \mid (\exists v' \in \ScatFact_1(u') : v' \in \ScatFact(w'))\land p(u)\leq p(w')\}\\
		\iff& \{u \in \Sigma^k \mid (\exists v\in\Sigma:\, v\in \al(u) \cap \al(w))\land p(u)\leq p(w)\} \\
		&= \{u' \in \Sigma^k \mid (\exists v'\in\Sigma:\,v' \in \al(u') \cap \al(w'))\land p(u')\leq p(w')\}\\
		\iff& p(w) = p(w').
	\end{align*}
	The last equivalence holds because we can choose $u=a_i^{p(w)[i]}$ and $u'=a_i^{p(w')[i]}$ for all $a_i\in\al(w)$. \qed
\end{proof}

%% file: universality.tex
In this section, we want to investigate Problem~\ref{problem:univ}, i.e., reach $k$-universality w.r.t. some word with as few
jumbles as possible. Recall that $\ell<|w|$ and $\iota(w)<k\leq|w|$.
For instance, for $w = \ta^2\tb^2\tc^2$, we have $\SJSF_2(w, 0) \cup \SJSF_2(w, 1)=\Sigma^2$.
Note that we can only get $\Sigma^k$ if all letters occur at least $k$ times in $w$, i.e., $k\in[\iota(w),\min_{\ta\in\Sigma}|w|_{\ta}]$. Set $I=[\iota(w),\min_{\ta\in\Sigma}|w|_{\ta}]$ for abbreviation.

\ifpaper
\input{grow}
\input{PROOFgrow}
\else

\input{grow}
\input{PROOFgrow}
\fi

With \Cref{grow} we can show that we need at most
$k - \iota(w)$ many jumbles.

\ifpaper
\input{lowerl}
\input{PROOFlowerl}

\else

\input{lowerl}
\input{PROOFlowerl}
\fi

The next lemma shows that in each iteration, we get new scattered factors, if the jumbled universality is not reached yet.

\ifpaper
\input{onemore}
\input{PROOFonemore}
\else

\input{onemore}
\input{PROOFonemore}
\fi

Note that $k-\iota(w)$ is not a tight bound for $\ell$ witnessed by $w'=\ta^3\tb^3$ resulting in $\SJSF_3(w,0)\cup\SJSF_3(w,1)=\Sigma^3$. The following part connects the classical $k$-universality with the jumbled $k$-universality. Similarly, as in \Cref{sec:jumblesimon} we observe that the letters that occur more than $k$ times play an important role in the following results.

\ifpaper
\input{onejumble}
\input{PROOFonejumble}
\else

\input{onejumble}
\input{PROOFonejumble}
\fi

\ifpaper
\input{jumbleiota+1}
\else
\input{jumbleiota+1}
 \input{PROOFjumbleiota+1}
\fi


\ifpaper
\input{wx}

\input{PROOFwx}
\else

\input{wx}
\input{PROOFwx}
\fi

Note that $w_x$ is not unique for a given $x\in\cmpl{\al(r(w))}$: given $\mathtt{aabc.bcbc}$
we may move the $\ta$ in front of and at the end of $r(w)$. If $x$ is the only absent letter in the rest, then putting it in front of the rest yields a shortest new arch. For \Cref{wx} it suffices to pick one. Nevertheless, we believe that it is worth to apply this approach of moving letters on one such $w_x$ for a given $x$. Before we define a refined factorisation, we give an example.
Consider the word $\mathtt{a^3bcd.ab^3cd.cd^2cd}$ over $\Sigma=\{\ta,\tb,\tc,\td\}$. The alphabet of the rest (in the first iteration) is $\Sigma_1=\{\tc,\td\}$ and we factorise the rest w.r.t. $\al(r(w))$ obtaining $r(w)=\tc\td.\td\tc.\td$ with $r^1(w)=\td$. Thus, $\tc\td$ is the first arch of $\r(w)$  which we can fill up with {\em superfluous} letters
from $w$'s arches. Now, we iterate through $\Sigma$: considering $\ta$, we notice that in $\ar_i(w)$, we have three $\ta$ and by \Cref{wx} we can move them towards the end in order to increase the universality. Thus here, we get $\mathtt{abcd.ab^3cd.acd.adc.d}$. For $\tb\in\Sigma$, we apply the same idea obtaining $\mathtt{abcd.abcd.bacd.badc.d}$ - a $4$-universal word. Thus, $w$ itself is jumbled $4$-universal.


For easier readability, introduce first the definition of the {\em potential} using the observation that only $k\in[\iota(w),\min_{\ta\in\Sigma}|w|_{\ta}]$ needs to be considered.

\begin{definition}
	Define the potential of $w\in\Sigma^{\ast}$ by $\nabla(w) = \min_{\ta \in \Sigma} \abs{w}_{\ta} - \iota(w)$.
\end{definition}

By the definition of $k$-universality, we have immediately $\nabla(w)\geq 0$. If $\al(w) \subset \Sigma$ then we have $\nabla(u) = 0$.

\ifpaper
\input{addedpotential}
\input{PROOFaddedpotential}
\else

\input{addedpotential}
\input{PROOFaddedpotential}
\fi

In this section, we started to investigate the jumbled scattered factor universality. These results are only a beginning of the research in this area.

%% file: grow.tex
\begin{lemma}\label{grow}
	Let $w \in \Sigma^*$ and $k + 1 \in I$. Then, we have for all $\ell \in [\abs{w}]$
	$$\JScatFact_k(w, \ell) = \Sigma^k \Rightarrow \JScatFact_{k+1}(w, \ell + 1) = \Sigma^{k + 1}.$$
\end{lemma}

%% file: PROOFgrow.tex
\begin{proof}
	By definition, we have $\JSF_{k+1}(w, \ell + 1) \subseteq \Sigma^{k+1}$. Let
	$u \in \Sigma^{k + 1}$, $v \in \Sigma^k$ and $x \in \Sigma$ with $u = vx$.
	Since $v \in \JSF_k(w, \ell)$ and $p(u) \leq p(w)$, we can apply \Cref{lemma:append} and get $\delta_w(u) \leq \delta_w(v) + 1 \leq \ell + 1$.
	By \Cref{lem:jsf_inclusion} we get $u \in \JSF_{k+1}(w, \ell + 1)$.\qed
\end{proof}

%% file: lowerl.tex
\begin{proposition}\label{lowerl}
	For all $w \in \Sigma^*$ there exist $\ell\leq k-\iota(w)$ and 
	$k \in I$, such that
	$\bigcup_{\ell' \in [0, \ell]} \JScatFact_k(w, \ell') = \Sigma^k$.
\end{proposition}

%% file: PROOFlowerl.tex
\begin{proof}
	We only need to consider $\ell = k - \iota(w)$ since for any larger $\ell$
	the union of jumbled scattered factors just becomes larger as well. We show this claim by
	induction over $k$. If we consider $k = \iota(w)$, then we can choose $\ell = 0$ and get
	$\JScatFact_k(w, 0) = \ScatFact_k(w) = \Sigma^k$
	since $w$ is $k$-universal. 
	
	Now assume
	$\bigcup_{\ell' \in [0, \ell]} \JScatFact_k(w, \ell') = \Sigma^k \quad \text{with } \ell = k - \iota(w)$
	for some $k \in [\iota(w), \min_{\ta \in \Sigma} \abs{w}_{\ta} - 1]$. By \Cref{grow}, we can conclude
	$$\bigcup_{\ell' \in [0, \ell + 1]} \JScatFact_{k + 1}(w, \ell') = \Sigma^{k + 1} \quad \text{with } \ell + 1 = k + 1 - \iota(w).\quad\qed$$
\end{proof}

%% file: onemore.tex
\begin{proposition}\label{onemore}
	Let $w \in \Sigma^*$, $k \in I$ and  $\JSF_k(w, \ell) \subset \Sigma^k$. Then, we have
	$\JSF_k(w, \ell + 1) \setminus \JSF_k(w, \ell) \neq \emptyset.$
\end{proposition}

%% file: PROOFonemore.tex
\begin{proof}
	We show this claim by {\em dejumbling} an arbitrary word that needs more than
	$\ell$ jumbles.
	Let $u \in \Sigma^k \setminus \JSF_k(w, \ell)$. This implies $\delta_w(u) > \ell$. If
	$\delta_w(u) = \ell + 1$, then we are done. Thus, we assume
	$\delta_w(u) > \ell + 1$. Let $x, y \in \Sigma^*, \ta \in \Sigma$ such that
	$u = x\ta y$ and $\delta_w(xy) = \delta_w(u) - 1$. Furthermore, let
	$\alpha, \beta \in \Sigma^*$ such that $w = \alpha \beta$ and
	$\lcsf(x, \alpha) + \lcsf(y, \beta) = \lcsf(u, w)$. W.l.o.g. we assume
	$\abs{\alpha}_\ta > \abs{x}_\ta$ (otherwise consider $\beta$ and $y$). Then,
	there exist $\mu, \nu \in \Sigma^*$
	such that $x = \mu \nu$ and $\lcsf(\alpha, \mu \ta \nu) = \lcsf(\alpha, x)+1$.
	Now, we have $u' = \mu \ta \nu y \notin \JSF_k(w, \ell)$ with $\abs{u'} = k$ and
	$\delta_w(u') = \delta(u) - 1$. We can repeat this process $\delta_w(u) - (\ell + 1)$
	times until we get a word that is $(\ell + 1)$-jumbled.\qed
\end{proof}

%% file: onejumble.tex
\begin{lemma}\label{onejumble}
  For $w \in \Sigma^*$ with $\iota(w)=k$ and $|w|_{\ta} \geq k + 1$ for all $\ta\in\Sigma$, we have
  $\SJScatFact_{k + 1}(w, 1) = (\ScatFact_k(w) \cdot \cmpl{\al(r(w))}) \setminus \ScatFact_{k+1}(w).$
\end{lemma}

%% file: PROOFonejumble.tex
\begin{proof}
	Let $u \in \SJScatFact_{k+1}(w, 1)$. Then, $u[1..k] \in \ScatFact_k(w)$ since
	$w$ is $k$-universal. If $u[k+1]$ were in $\al(r(w))$, then
	$\delta_w(u) = 0 \neq 1$.
	
	\medskip
	
	Let $u \in \ScatFact_k(w) \cdot \cmpl{\al(r(w))}$ and $u \notin \ScatFact_{k+1}(w)$.
	Then, $u$ must have at least one jumble and it cannot have more than one, since
	$u[1..k] \in \ScatFact_k(w)$. Furthermore, $$\norm{p(u)}_\infty \leq \norm{p(u[1..k])}_\infty + \norm{p(u[k+1])}_\infty = k + 1.$$
	By $|w|_{\ta}\geq k+1$ for all $\ta\in\Sigma$, each entry of $p(w)$ is at least $k+1$.
	This concludes the proof.  \qed
\end{proof}

%% file: jumbleiota+1.tex
\begin{proposition}\label{jumbleiota+1}
  For $w \in \Sigma^*$ with $\iota(w)=k$ and $|w|_{\ta} \geq k + 1$ for all $\ta\in\Sigma$, we have $\SJScatFact_{k + 1}(w, 0) \cup \SJScatFact_{k + 1}(w, 1) = \Sigma^{k + 1}$.
\end{proposition}

%% file: PROOFjumbleiota+1.tex
\begin{proof}
  The following equation holds
  \begin{align*}
    & \SJScatFact_{k + 1}(w, 0) \cup \SJScatFact_{k + 1}(w, 1) \\
    =& \ScatFact_{k + 1}(w) \cup (\ScatFact_k(w) \cdot \cmpl{\al(r(w))}) \setminus \ScatFact_{k+1}(w) \\
    =& \ScatFact_{k + 1}(w) \cup \ScatFact_k(w) \cdot \cmpl{\al(r(w))} \\
    =& \ScatFact_{k}(w) \cdot \al(r(w)) \cup \ScatFact_k(w) \cdot \cmpl{\al(r(w))} \\
    =& \Sigma^{k+1}.\quad\qed
  \end{align*}
\end{proof}

%% file: wx.tex
\begin{lemma}\label{wx}
  Let $w \in \Sigma^*$ with $\iota(w)=k$ and $|w|_{\ta}\geq k+1$ for all $\ta\in\Sigma$. For all $x \in \cmpl{\al(r(w))}$ there exists $w_x \in \JScatFact(w, 1)$ such that $\ScatFact_{k+1}(w) \cup \bigcup_{x \in \cmpl{\al(r(w))}} \ScatFact_{k+1}(w_x) = \Sigma^{k+1}$.
\end{lemma}

%% file: PROOFwx.tex
\begin{proof}
	Let $w \in \Sigma^*$ be $k$-universal and $|w|_{\ta} \geq k + 1$ for all $\ta\in\Sigma$.
	Let $x \in \cmpl{\al(r(w))}$. Since $|w|_{\ta}\geq k+1$ and $w$ is just
	$k$-universal, there exists a smallest $i \in [k]$ with $\abs{\ar_i(w)}_x > 1$.
	Let $\alpha, \beta \in \Sigma^*$ with $x \notin \al(\alpha)$ such that $\ar_i(w) = \alpha x\beta$. We
	define $$w_x = \ar_1(w) \cdots \ar_{i-1}(w) \cdot \alpha \cdot \beta \cdot \ar_{i+1}(w) \cdots \ar_k(w) \cdot x \cdot r(w).$$
	We can write $w = uxv \cdot r(w)$ and $w_x = uvx \cdot r(w)$ for some
	$u, v \in \Sigma^*$. Then, we have
	\[
	\lcsf(w, w_x) = \lcsf(uxv \cdot r(w), uvx \cdot r(w)) = \abs{u} + \lcsf(xv, vx) + \abs{r(w)}.
	\]
	Note that $\lcsf(xv,vx)\geq|v|$.\\
	\textbf{Supposition:} $\lcsf(xv,vx)=|xv|$\\
	This is equivalent to $xv=vx$. By the Lyndon-Sch{\" u}tzenberger Theorem, also known as xy-yx-lemma, we get $v=x^{|v|}$. This implies $i=k$ and $v=\beta$. Thus, $|\ar_{k}(w)|_{x}=2$ and $m(w)[\abs{m(w)}] = x$ - a contradiction
	to the arch factorisation.\\
	Hence we have $\lcsf(xv,vx)=|v|$ and consequently 
	\[
	\lcsf(w,w_x)=|u|+|v|+|r(w)|=|w|-1.
	\]
	This implies $\delta_w(w_x)=1$ as claimed.
	
	\medskip
	
	For the second claim the $\subseteq$-inclusion is immediate.  Let $u \in \Sigma^{k+1}$. If $u$ is a scattered factor $w$,
	then the statement trivially holds. Thus, we consider the case that $u$ is
	not a scattered factor of $w$. Let $v \in \Sigma^*$ and $x \in \Sigma$ such
	that $u = vx$. Then $u \in \ScatFact_{k+1}(w_x)$ since either $\iota(w_x) = k + 1$, or
	$\iota(w_x) = k$ and $x \in \al(r(w_x))$ by construction.\qed
\end{proof}

%% file: addedpotential.tex
\begin{lemma}
	Let $u, v \in \Sigma^*$ with $\min_{\ta \in \Sigma} \abs{uv}_\ta = \min_{\ta \in \Sigma} \abs{u}_\ta + \min_{\ta \in \Sigma} \abs{v}_\ta$. Then,  $\nabla(u) + \nabla(v)\in\{\nabla(uv),\nabla(uv)-1\}$.
\end{lemma}

%% file: PROOFaddedpotential.tex
\begin{proof}
	We have that there exists $\ta\in\Sigma$ with $\min_{x \in \Sigma} \abs{u}_x =|u|_\ta$ and $\min_{x \in \Sigma} \abs{v}_x=|v|_\ta$. Thus, we have
	\begin{align*}
		\nabla(u)+\nabla(v)=|u|_\ta-\iota(u)+|v|_\ta-\iota(v)=|uv|_\ta-(\iota(u)+\iota(v)).
	\end{align*}
	Since $\ta$ occurs minimally often in $u$ and in $v$ resp. it occurs minimally often in $uv$.
	The claim follows by \cite{fleischmann2021scattered}.\qed
\end{proof}

%% file: conclusion.tex
In this work, we introduced jumbled scattered factors as a combination of
(absent) scattered factors and jumbled words. After fundamental properties on
the jumbled index, we continued the investigation of the jumbled Simon's
equivalence and the jumbled $k$-universality based on the respective versions of
scattered factors. For the equivalence relation, it is left as future work to give a more fine-grained analysis of the classes that contain words with the same Parikh vector and to
determine its number of equivalence classes which might lead to insights of the
number of congruence classes for the original Simon's congruence.
Regarding the universality, we are also interested in the following problem: 
Given $w \in \Sigma^*$ with $\ell \in [|w|]$. What
is the maximal value $k$ s.t. $\bigcup_{\ell' \in [0,\ell]}\JScatFact_{k}(w,\ell') = \Sigma^k$?
In light of the recent interest of universality from a language perspective this
and Problem~\ref{problem:univ} could be extended to the inspection of not only
a single but a set of words.

Further, one can extend classical scattered factor problems from combinatorics
on words to their jumbled equivalent, e.g., finding the longest common jumbled
scattered factor for a given pair of words. 

%